\documentclass[12pt]{amsart}

\usepackage{a4wide}

\usepackage{amssymb,amsmath,stmaryrd,xy,graphics}

\xyoption{all}

\newtheorem{theorem}{Theorem}[section]

\theoremstyle{remark}

\begin{document}
\title{Hyperpolar Homogeneous Foliations on\\
Symmetric Spaces of Noncompact Type}

\author{J\"{u}rgen Berndt}

\begin{abstract}
A foliation ${\mathcal F}$ on a Riemannian manifold $M$ is homogeneous
if its leaves coincide with the orbits of an isometric action on $M$.
A foliation ${\mathcal F}$ is polar if it admits a section, that is, a connected closed
totally geodesic submanifold of $M$ which intersects each leaf of ${\mathcal F}$, and intersects
orthogonally at each point of intersection.
A foliation ${\mathcal F}$ is hyperpolar if it admits a flat section.
These notes are related to joint work with Jos\'{e} Carlos D\'{i}az-Ramos
and Hiroshi Tamaru about hyperpolar homogeneous foliations on Riemannian
symmetric spaces of noncompact type.
Apart from the classification result which we proved in \cite{BDT}, we
present here in more detail some relevant material about symmetric spaces
of noncompact type,
and discuss the classification in more detail for the special case
$M = SL_{r+1}({\mathbb R})/SO_{r+1}$.
\end{abstract}


\maketitle


\section{Isometric actions on the real hyperbolic plane}

The special linear group
$$
G = SL_2({\mathbb R}) =
\left\{ \left.\begin{pmatrix} a & b \\ c & d \end{pmatrix} \right|  a,b,c,d \in {\mathbb R},\ ad-bc = 1 \right\}
$$
acts on the upper half plane
$$ H = \{ z \in {\mathbb C} \mid \Im(z) > 0 \}$$
by linear fractional transformations of the form
$$ \begin{pmatrix} a & b \\ c & d \end{pmatrix} \cdot z = \frac{az+b}{cz+d}.$$
By equipping $H$ with the Riemannian metric
$$
ds^2 = \frac{dx^2 + dy^2}{y^2} = \frac{dzd\bar{z}}{\Im(z)^2}\ \ (z = x+iy)
$$
we obtain the Poincar\'{e} half-plane model of the real hyperbolic plane ${\mathbb R}H^2$.
The above action of $SL_2({\mathbb R})$ is isometric with respect to this metric.
The action is transitive and the stabilizer at $i$ is
$$
K = \left\{ \left.\begin{pmatrix} \cos(s) & \sin(s) \\ -\sin(s) & \cos(s) \end{pmatrix} \right|  s \in {\mathbb R} \right\} = SO_2.
$$
The subgroup $K$ is a maximal compact subgroup of $G$, and
 we can identify ${\mathbb R}H^2$ with the homogeneous space $G/K = SL_2({\mathbb R})/SO_2$ in the usual way.
There is a unique decomposition of every matrix in $SL_2({\mathbb R})$ into three matrices of the form
$$
\begin{pmatrix} a & b \\ c & d \end{pmatrix}
= \begin{pmatrix} \cos(s) & \sin(s) \\ -\sin(s) & \cos(s) \end{pmatrix}
\begin{pmatrix} \exp(t) & 0 \\ 0 & \exp(-t) \end{pmatrix}
\begin{pmatrix} 1 & u \\ 0 & 1 \end{pmatrix}\ \ (s,t,u \in {\mathbb R}).
$$
We define an abelian subgroup $A$ of $G$ by
$$
A = \left\{ \left.\begin{pmatrix} \exp(t) & 0 \\ 0 & \exp(-t) \end{pmatrix} \right|  t \in {\mathbb R} \right\}
$$
and a nilpotent subgroup $N$ of $G$ by
$$
N = \left\{ \left.\begin{pmatrix} 1 & u \\ 0 & 1 \end{pmatrix} \right|  u \in {\mathbb R} \right\}.
$$
Then $G$ is diffeomorphic to the manifold product $KAN$. The decomposition $G = KAN$ is known as
an Iwasawa decomposition of $G$. Each of three subgroups $K$, $A$ and $N$ acts isometrically on
${\mathbb R}H^2$.

\begin{itemize}
\item The orbits of $K$ consist of the single point $i$, which is the fixed point of the action of $K$ on ${\mathbb R}H^2$,
and the geodesic circles in ${\mathbb R}H^2$ whose center is at $i$.
\item The orbits of $A$ are the geodesic ${\mathbb R} \to {\mathbb R}H^2,t \mapsto \exp(2t)i$ and the
curves in ${\mathbb R}H^2$ which are equidistant to this geodesic.
\item The orbits of $N$ are the horospheres given by the equations $\Im(z) = y$, $y \in {\mathbb R}_+$.
\end{itemize}

The action of $K$ has exactly one singular orbit, and the orbits of $A$ and $N$ form a Riemannian foliation
on the real hyperbolic plane. Every isometric action on ${\mathbb R}H^2$ by a nontrivial connected subgroup $H$
of $G$ is orbit equivalent to one of these three actions, that is, there is an isometry of ${\mathbb R}H^2$
which maps the orbits of the action of $H$ onto the orbits of the action of either $K$, or of $A$, or of $N$.

The real hyperbolic plane is the most elementary example of a Riemannian symmetric space
of noncompact type. A joint project with Jos\'{e} Carlos D\'{i}az Ramos
(Santiago de Compostela) and Hiroshi Tamaru (Hiroshima)
is to generalize aspects of what we described above from ${\mathbb R}H^2$
to Riemannian symmetric spaces of noncompact type.

\section{The problem}

There are several natural ways to generalize what we discussed in the previous section. The cohomogeneity
of each of the three actions is equal to one. The cohomogeneity of an isometric action is the codimension of
a principal orbit of the action, or equivalently, the dimension of the orbit space of the action.
A cohomogeneity one action is therefore an action for which the principal orbits of the action are hypersurfaces.

Cohomogeneity one actions are special cases of polar and hyperpolar actions.
Let $M$ be a connected complete Riemannian manifold and $H$ a
connected closed subgroup of the isometry group $I(M)$ of $M$. Then
each orbit $H \cdot p = \{ h(p) \mid h \in H\}$, $p \in M$, is a
connected closed submanifold of $M$. A connected complete
submanifold ${\mathcal S}$ of $M$ that meets each orbit of the
$H$-action and intersects the orbit $H \cdot p$ perpendicularly at
each point $p \in {\mathcal S}$ is called a section of the action. A
section ${\mathcal S}$ is always a totally geodesic submanifold of
$M$ (see e.g.\ \cite{HLO}). In general, actions do not admit a
section. The action of $H$ on $M$ is called polar if it has a
section, and it is called hyperpolar if it has a flat section. For
motivation and classification of polar and hyperpolar actions on
Euclidean spaces and symmetric spaces of compact type we refer to
the papers by Dadok \cite{Da85}, Podest\`{a} and Thorbergsson
\cite{PT99}, and Kollross \cite{Ko02}, \cite{Ko07}. If all orbits of
$H$ are principal, then the orbits form a homogeneous foliation
${\mathcal F}$ on $M$. In general, a foliation ${\mathcal F}$ on $M$
is called homogeneous if the subgroup of $I(M)$ consisting of all
isometries preserving ${\mathcal F}$ acts transitively on each leaf
of ${\mathcal F}$. Homogeneous foliations are basic examples of
metric foliations. A homogeneous foliation is called polar resp.\
hyperpolar if its leaves coincide with the orbits of a polar resp.\
hyperpolar action.

The problem I want to discuss in this article is the classification of
hyperpolar homogeneous foliations on Riemannian symmetric spaces of noncompact type.

\section{Riemannian symmetric spaces of noncompact type}

This section contains a brief introduction to Riemannian symmetric spaces
of noncompact type. A basic reference is the book \cite{H78} by Sigurdur Helgason.

A Riemannian symmetric space is a connected Riemannian manifold $M$ with the property that
for each point $p \in M$ there exists an isometric involution $s_p \in I(M)$
for which $p$ is an isolated fixed point. Geometrically, such an involution $s_p$ is
the geodesic reflection of $M$ in the point $p$. Every Riemannian symmetric space is
homogeneous, that is, the isometry group $I(M)$ acts transitively on $M$. The Riemannian
universal covering space $\widetilde{M}$ of a
Riemannian symmetric space $M$ is also a Riemannian symmetric space, and
each factor in the de Rham decomposition
$$
\widetilde{M} = M_0 \times M_1 \times \cdots \times M_k
$$
is a Riemannian symmetric space as well. Here $M_0$ is the Euclidean factor (which may be $0$-dimensional),
and each manifold $M_i$, $1 \leq i \leq k$, is an irreducible simply connected Riemannian symmetric space.
A Riemannian symmetric space $M$ is said to be of compact type if $\dim M_0 = 0$ and each $M_i$ is compact,
and it is said to be of noncompact type if $\dim M_0 = 0$ and each $M_i$ is noncompact.
Every Riemannian symmetric space of noncompact type is simply connected.
The rank of a Riemannian symmetric space is the maximal dimension of a flat totally geodesic
submanifold.

Riemannian symmetric spaces have been classified by \'{E}lie Cartan.
The Riemannian symmetric spaces of noncompact type
and with rank one are the hyperbolic spaces ${\mathbb F}H^n$ over the normed real division
algebras ${\mathbb F} \in \{{\mathbb R},{\mathbb C},{\mathbb H},{\mathbb O}\}$, where $n \geq 2$ and $n = 2$ if
${\mathbb F} = {\mathbb O}$.

Let $M$ be a connected Riemannian symmetric space of noncompact
type. We denote by $n$ the dimension of $M$ and by $r$ the rank of
$M$. It is well known that $M$ is diffeomorphic to ${\mathbb R}^n$.
Let $G = I^o(M)$ be the connected component of $I(M)$
containing the identity transformation of $M$. Then $G$ is a connected
semisimple real Lie group with trivial center.
Let $o \in M$ and $K = \{g \in G \mid g(o) = o\}$ be the isotropy subgroup of $G$ at $o$.
Then $K$ is a maximal compact subgroup of $G$.
We identify $M$ with the homogeneous space $G/K$ in the
usual way and denote by ${\mathfrak g}$ and ${\mathfrak k}$ the Lie algebra of $G$
and $K$, respectively. Let
$$
B : {\mathfrak g} \times {\mathfrak g} \to {\mathfrak g}\ ,\ (X,Y) \mapsto {\rm tr}({\rm ad}(X) \circ
{\rm ad}(Y))
$$
be the Killing form of ${\mathfrak g}$. It is well-known from Cartan's criterion for semisimple
Lie algebras that the Killing form $B$ is non-degenerate. Let
$${\mathfrak p} = \{ X \in {\mathfrak g} \mid B(X,Y) = 0\ {\rm for\ all}\ Y \in {\mathfrak k}\}
$$
be the orthogonal complement of ${\mathfrak k}$ with respect
to $B$. Then $${\mathfrak g} = {\mathfrak k} \oplus {\mathfrak p}$$ is a Cartan decomposition of
${\mathfrak g}$, that is, we have
$$
[{\mathfrak k},{\mathfrak k}] \subset {\mathfrak k}\ ,\ [{\mathfrak k},{\mathfrak p}] \subset {\mathfrak p}\ ,\
[{\mathfrak p},{\mathfrak p}] \subset {\mathfrak k},
$$
and $B$ is negative definite on ${\mathfrak k}$ and positive definite on ${\mathfrak p}$.
The Cartan involution corresponding to the Cartan decomposition ${\mathfrak g} =
{\mathfrak k} \oplus {\mathfrak p}$ is the automorphism of ${\mathfrak g}$ defined by
$\theta(X) = X$ for all $X \in {\mathfrak k}$ and $\theta(X) = -X$ for all $X \in {\mathfrak p}$.
We can
define a positive definite inner product on ${\mathfrak g}$ by $\langle
X,Y \rangle = -B(X,\theta Y)$ for all $X,Y \in {\mathfrak g}$. As usual, we identify
${\mathfrak p}$ with $T_oM$, and we normalize the Riemannian metric on $M$ so
that its restriction to $T_o M\times T_o M={\mathfrak p}\times{\mathfrak p}$
coincides with $\langle\,\cdot\,,\,\cdot\,\rangle$.

{\bf Example.} $M = SL_{r+1}({\mathbb R})/SO_{r+1}$. The linear group $GL_{r+1}({\mathbb R})$
can be identified with the group of all regular $(r+1) \times (r+1)$-matrices with real
coefficients, or equivalently, with the group of all automorphisms of ${\mathbb R}^{r+1}$.
The special linear group $G = SL_{r+1}({\mathbb R})$ of all orientation- and volume-preserving transformations
in $GL_{r+1}({\mathbb R})$ is
$$ SL_{r+1}({\mathbb R}) = \{A \in GL_{r+1}({\mathbb R}) \mid \det(A) = 1\}.$$
The orthogonal group $O_{r+1}$ is the group of transformation on ${\mathbb R}^{r+1}$ preserving the
standard inner product
$$ \langle x , y \rangle = \sum_{i=1}^{r+1} x_iy_i $$
on ${\mathbb R}^{r+1}$. The special orthogonal group $K = SO_{r+1} = O_{r+1} \cap SL_{r+1}({\mathbb R})$ is the subgroup of $O_{r+1}$
preserving the standard orientation of ${\mathbb R}^{r+1}$. It is the connected component of $O_{r+1}$ containing the identity
transformation of ${\mathbb R}^{r+1}$ and can be identified with
$$SO_{r+1} = \{A \in SL_{r+1}({\mathbb R}) \mid A^t = A^{-1}\},$$
where $A^t$ is the transpose of $A$ and $A^{-1}$ is the inverse of $A$.
The subgroup $SO_{r+1}$ is a maximal compact subgroup of $SL_{r+1}({\mathbb R})$.

The Lie algebra ${\mathfrak g}{\mathfrak l}_{r+1}({\mathbb R})$ of $GL_{r+1}({\mathbb R})$ is
the $({r+1})^2$-dimensional real vector space of all $({r+1}) \times ({r+1})$-matrices with real coefficients
together with the Lie algebra structure defined by
$$ [X,Y] = XY - YX\ ,\ X,Y \in {\mathfrak g}{\mathfrak l}_{r+1}({\mathbb R}).$$
The Lie algebra ${\mathfrak g} = {\mathfrak s}{\mathfrak l}_{r+1}({\mathbb R})$ is the subalgebra of
${\mathfrak g}{\mathfrak l}_{r+1}({\mathbb R})$ defined by
$$ {\mathfrak s}{\mathfrak l}_{r+1}({\mathbb R}) = \{X \in {\mathfrak g}{\mathfrak l}_{r+1}({\mathbb R}) \mid
{\rm trace}(X) = 0\},$$
and the Lie algebra ${\mathfrak k} = {\mathfrak s}{\mathfrak o}_{r+1}$ is
$$ {\mathfrak s}{\mathfrak o}_{r+1} = \{ X \in {\mathfrak s}{\mathfrak l}_{r+1}({\mathbb R}) \mid X^t = -X\}.$$

The orthogonal complement ${\mathfrak p}$ of ${\mathfrak k} = {\mathfrak s}{\mathfrak o}_{r+1}$
in ${\mathfrak g} = {\mathfrak s}{\mathfrak l}_{r+1}({\mathbb R})$ with
respect to the Killing form $B$ of ${\mathfrak s}{\mathfrak l}_{r+1}({\mathbb R})$
is the real vector space of all symmetric $({r+1}) \times ({r+1})$-matrices, that is,
$$ {\mathfrak p} = \{ X \in {\mathfrak s}{\mathfrak l}_{r+1}({\mathbb R}) \mid X^t = X\}.$$
The Cartan decomposition ${\mathfrak g} = {\mathfrak k} \oplus {\mathfrak p}$ therefore
corresponds to the decomposition of $X \in {\mathfrak s}{\mathfrak l}_{r+1}({\mathbb R})$ into
its skewsymmetric and its symmetric part.

The Cartan involution corresponding to the above Cartan
decomposition of ${\mathfrak s}{\mathfrak l}_{r+1}({\mathbb R})$ is
$$ \theta : {\mathfrak s}{\mathfrak l}_{r+1}({\mathbb R}) \to {\mathfrak s}{\mathfrak l}_{r+1}({\mathbb R})\ ,\ X \mapsto -X^t.$$
It is obvious that ${\mathfrak s}{\mathfrak o}_{r+1}$ is the fixed point set
of $\theta$.

\section{Lie triple systems and totally geodesic submanifolds}

A linear subspace ${\mathfrak m}$ of ${\mathfrak p}$ is called a Lie triple
system if
$$
[[{\mathfrak m},{\mathfrak m}],{\mathfrak m}] \subset {\mathfrak m}.
$$
If ${\mathfrak m} \subset {\mathfrak p}$ is a Lie triple system, then
$$
{\mathfrak g}_{\mathfrak m} = [{\mathfrak m},{\mathfrak m}] \oplus {\mathfrak m} \subset
{\mathfrak k} \oplus {\mathfrak p} = {\mathfrak g}
$$
is a subalgebra of ${\mathfrak g}$. Let $G_{\mathfrak m}$ be the connected subgroup of
$G$ with Lie algebra ${\mathfrak g}_{\mathfrak m}$. Then the orbit
$$
F_{\mathfrak m} = G_{\mathfrak m} \cdot o
$$
of $G_{\mathfrak m}$ through $o$ is a connected complete totally geodesic submanifold of $M$.

Conversely, let $F$ be a totally geodesic submanifold of $M$ with $o \in F$.
Then the tangent space $T_oF$ of $F$ at $o$ can be considered as a linear subspace
of ${\mathfrak p}$ through the identification $T_oM \cong {\mathfrak p}$, and $T_oF$
is a Lie triple system.

Therefore there exists a one-to-one correspondence between Lie triple systems in ${\mathfrak p}$
and connected complete totally geodesic submanifolds of $M$ containing the point $o$.

\section{The restricted root space decomposition}

Let ${\mathfrak a}$ be a maximal abelian subspace of ${\mathfrak p}$ and denote
by ${\mathfrak a}^*$ the dual space of ${\mathfrak a}$. For each $\lambda\in {\mathfrak a}^*$
we define
$${\mathfrak g}_\lambda=\{X\in {\mathfrak g} \mid {\rm ad}(H)X=\lambda(H)X \mbox{ for
all }H\in {\mathfrak a}\}.$$
A restricted root of ${\mathfrak g}$ is an element $0 \neq \lambda \in {\mathfrak a}^*$
such that ${\mathfrak g}_\lambda\neq \{0\}$. We denote by
$\Sigma$ the set of all restricted roots of ${\mathfrak g}$.
Since ${\mathfrak a}$ is abelian,
${\rm ad}({\mathfrak a})$ is a commuting family of selfadjoint linear
transformations of ${\mathfrak g}$.
The corresponding eigenspace decomposition
$${\mathfrak g}= {\mathfrak g}_0\oplus\left(\bigoplus_{\lambda\in\Sigma}{\mathfrak g}_\lambda\right)$$
is an orthogonal direct sum and called the restricted root space
decomposition of ${\mathfrak g}$ determined by ${\mathfrak a}$. The eigenspace
${\mathfrak g}_0$ of $0 \in {\mathfrak a}^*$ is equal to
$${\mathfrak g}_0={\mathfrak k}_0\oplus{\mathfrak a},$$
where ${\mathfrak k}_0=Z_{{\mathfrak k}}({\mathfrak a})$ is
the centralizer of ${\mathfrak a}$ in ${\mathfrak k}$.

For each $\lambda\in {\mathfrak a}^*$
let $H_\lambda\in {\mathfrak a}$ denote the dual vector in ${\mathfrak a}$ with
respect to the inner product $\langle \cdot , \cdot \rangle$, that is, $\lambda(H)=\langle
H_\lambda,H\rangle$ for all $H\in {\mathfrak a}$. We then get an inner product
on ${\mathfrak a}^*$ by setting $\langle\lambda,\mu\rangle=\langle
H_\lambda,H_\mu\rangle$ for all $\lambda,\mu\in{\mathfrak a}^*$. Below we always consider ${\mathfrak a}^*$
to be equipped with this inner product.

{\bf Example.} $M = SL_{r+1}({\mathbb R})/SO_{r+1}$.
For $a = (a_1,\ldots,a_{r+1}) \in {\mathbb R}^{r+1}$ we denote by $\Delta(a) \in {\mathfrak g}{\mathfrak l}_{r+1,r+1}({\mathbb R})$ the
matrix with $\Delta(a)_{ii} = a_i$ and $\Delta(a)_{ij} = 0$
for all distinct $i,j \in \{1,\ldots,r+1\}$. Then
$$
{\mathfrak a}  =  \{ \Delta(a) \mid a = (a_1,\ldots,a_{r+1}) \in {\mathbb R}^{r+1},\ a_1 + \ldots + a_{r+1} = 0\}
$$
is a maximal abelian subspace of ${\mathfrak p}$.
Let $e_i \in {\mathfrak a}$ be the vector in ${\mathfrak a}$ which is obtained by putting
$a_i = 1$ and $a_j = 0$ for all $j \neq i$, $i,j \in \{1,\ldots,r+1\}$. Denote by $e_1^*,\ldots,e_{r+1}^* \in
{\mathfrak a}^*$ the dual vectors of $e_1,\ldots,e_{r+1}$. Then
$e_i^* - e_j^*$ for $i \neq j$ and $i,j \in \{1,\ldots,r+1\}$ is a restricted root of
${\mathfrak g} = {\mathfrak s}{\mathfrak l}_{r+1}({\mathbb R})$ with corresponding root space
$${\mathfrak g}_{e_i^* - e_j^*}  = \{ xE_{ij} \mid x \in {\mathbb R}\},$$
where $E_{ij} \in {\mathfrak g}{\mathfrak l}_{r+1,r+1}({\mathbb R})$ is the matrix with
$(E_{ij})_{ij} = 1$ and zero everywhere else. We thus have
$$\Sigma = \{ e_i^* - e_j^* \mid i \neq j,\ i,j \in \{1,\ldots,r+1\} \}.$$
The centralizer of ${\mathfrak a}$ in ${\mathfrak k}$ is trivial, and therefore we have
$$ {\mathfrak g}_0 = {\mathfrak a}.$$
The restricted root space decomposition of ${\mathfrak s}{\mathfrak l}_{r+1,r+1}({\mathbb R})$
is therefore given by
$$
{\mathfrak s}{\mathfrak l}_{r+1,r+1}({\mathbb R}) = {\mathfrak a} \oplus \left(
\bigoplus_{\substack{i,j = 1 \\ i \neq j}}^{r+1} {\mathbb R}E_{ij}\right)
$$

\section{The Iwasawa decomposition}

A subset $\Lambda=\{\alpha_1,\dots,\alpha_r\} \subset \Sigma$
is called a set of simple roots of $\Sigma$ if every $\lambda \in \Sigma$
can be written in the form
$$ \lambda = \sum_{i=1}^r c_i\alpha_i $$
with some integers $c_1,\ldots,c_r \in {\mathbb Z}$ such that
$c_1,\ldots,c_r$ are either all nonpositive or all nonnegative.
A set of simple roots of $\Sigma$ always exists, and it is unique up to a transformation
in the Weyl group of $\Sigma$. The Weyl group of $\Sigma$ is the subgroup of
orthogonal transformations of ${\mathfrak a}^*$ which is generated by the reflections
$$
s_\lambda : {\mathfrak a}^* \to {\mathfrak a}^*\ ,\ X \mapsto X - 2
\frac{\langle X , \lambda \rangle}{|\lambda|^2}\lambda.
$$
The set
$$ \Sigma^+ = \{ \lambda \in \Sigma \mid \lambda = c_1\alpha_1 + \cdots + c_r\alpha_r,\
c_1,\ldots,c_r \geq 0\}$$
is called the set of positive restricted roots of $\Sigma$ with
respect to $\Lambda$. For the purpose of consistency in our article
our choice of simple roots will be the one used in \cite{K96}.

The subspace $${\mathfrak n}=\bigoplus_{\lambda\in\Sigma^+}{\mathfrak g}_\lambda$$ of
${\mathfrak g}$ is a nilpotent subalgebra of ${\mathfrak g}$. Moreover,
${\mathfrak a}\oplus{\mathfrak n}$ is a solvable subalgebra of ${\mathfrak g}$ and
${\mathfrak n}$ is the derived subalgebra of ${\mathfrak a} \oplus {\mathfrak n}$, that is,
$[{\mathfrak a} \oplus {\mathfrak n},{\mathfrak a} \oplus {\mathfrak n}] = {\mathfrak n}$. The
direct sum vector space decomposition
${\mathfrak g}={\mathfrak k}\oplus{\mathfrak a}\oplus{\mathfrak n}$ is called the Iwasawa
decomposition of ${\mathfrak g}$ with respect to ${\mathfrak a}$.
We emphasize that this is only a direct sum of vector spaces and not a direct
sum of Lie algebras.
Let $A$, $N$ and $AN$ be the connected
subgroups of $G$ with Lie algebra ${\mathfrak a}$, ${\mathfrak n}$ and
${\mathfrak a}\oplus {\mathfrak n}$, respectively. All three subgroups
$A$, $N$ and $AN$ are simply
connected and $G$ is diffeomorphic to the manifold product $K\times A\times
N$. Again, we emphasize that this is a product of smooth manifolds
and not a Lie group product.

It follows from the Iwasawa decomposition that the solvable Lie group $AN$ acts simply transitively
on $M$.  Therefore $M$ is isometric to the connected, simply connected
solvable Lie group $AN$ equipped with the left-invariant Riemannian
metric which is induced from the inner product
$\langle\,\cdot\,,\,\cdot\,\rangle$ restricted to ${\mathfrak a} \oplus {\mathfrak n}$.

{\bf Example.} $M = SL_{r+1}({\mathbb R})/SO_{r+1}$.
From the above description one can easily verify that
$$\alpha_i = e_i^* - e_{i+1}^*,\ i \in \{1,\ldots,r\},$$ form a set of simple roots for $\Sigma$.
We denote this set by $\Lambda$, that is, $\Lambda = \{\alpha_1,\ldots,\alpha_r\}$.
The corresponding set of positive restricted roots is
$$\Sigma^+  =  \{ e_i^* - e_j^* \mid i < j,\ i,j \in \{1,\ldots,r+1\} \}. $$
Using the above description of the restricted root spaces we see that ${\mathfrak n}$
is the nilpotent Lie algebra consisting of all strictly upper diagonal matrices, that is,
$${\mathfrak n} = \left\{ \begin{pmatrix}
0 & x_{12} & x_{13} & \cdots & x_{1,{r+1}} \\
0 & 0 & x_{23} & \cdots & x_{2,{r+1}} \\
\vdots & \vdots & \ddots & \cdots & \vdots \\
0 & 0 & 0 & \cdots & x_{r,{r+1}} \\
0 & 0 & 0 & \cdots & 0
\end{pmatrix} \right\}.$$
The Iwasawa decomposition of ${\mathfrak s}{\mathfrak l}_{r+1}({\mathbb R})$
therefore describes the unique decomposition of a matrix in
${\mathfrak s}{\mathfrak l}_{r+1}({\mathbb R})$ into the sum of a skewsymmetric matrix,
a diagonal matrix with trace zero, and a strictly upper triangular matrix.
On Lie group level, we get
$$N = \left\{ \begin{pmatrix}
1 & x_{12} & x_{13} & \cdots & x_{1,{r+1}} \\
0 & 1 & x_{23} & \cdots & x_{2,{r+1}} \\
\vdots & \vdots & \ddots & \cdots & \vdots \\
0 & 0 & 0 & \cdots & x_{r,{r+1}} \\
0 & 0 & 0 & \cdots & 1
\end{pmatrix} \right\}.$$
The Iwasawa decomposition of $SL_{r+1}({\mathbb R})$
therefore describes the unique decomposition of a matrix in
$SL_{r+1}({\mathbb R})$ into the product of an orthogonal matrix with determinant one,
a diagonal matrix with determinant one, and an upper triangular matrix with entries equal to one in
the diagonal.
The solvable Lie group $AN$ of upper triangular matrices with determinant one acts
simply transitively on the symmetric space $M = SL_{r+1}({\mathbb R})/SO_{r+1}$.

\section{The Dynkin diagram}

We now assign to the symmetric space $M$ a diagram consisting of
vertices, lines and arrows. Consider the set $\Lambda = \{\alpha_1,\ldots,\alpha_r\}$ of simple roots
of $\Sigma$. To each simple root $\alpha_i \in \Lambda$ we assign a vertex which we denote by
$ \xy  *\cir<2pt>{} \endxy $ if $2\alpha_i \notin \Sigma$ and by $ \xy  *\cir<2pt>{} *\cir<4pt>{} \endxy $ if
$2 \alpha_i \in \Sigma$.

One can show that the angle
between two simple roots in $\Lambda$ is one of the following four angles:
$$
\frac{\pi}{2}\ ,\ \frac{\pi}{3}\ ,\ \frac{\pi}{4}\ ,\ \frac{\pi}{6}.
$$
We connect the vertices corresponding to simple roots $\alpha_i$ and $\alpha_j$, $i \neq j$, by
$0$, $1$, $2$ or $3$ lines if the angle between $\alpha_i$ and $\alpha_j$ is
$\frac{\pi}{2}$, $\frac{\pi}{3}$, $\frac{\pi}{4}$ or $\frac{\pi}{6}$ respectively.
Moreover, if the vertices corresponding to $\alpha_i$ and $\alpha_j$ are connected by at least one line
and $\langle \alpha_i,\alpha_i \rangle > \langle \alpha_j,\alpha_j \rangle$, we draw an arrow
from the vertex $\alpha_i$ to the vertex $\alpha_j$. The resulting object is called
the Dynkin diagram associated with $M$.

We now list the Dynkin diagrams for the irreducible Riemannian symmetric spaces
of noncompact type . We view a root system as a subset of some Euclidean vector space $V$.
We also list for each symmetric space $M$ the multiplicities $(m_{\alpha_1},\ldots,m_{\alpha_r})$
of the simple roots $\alpha_1,\ldots,\alpha_r$.
\begin{itemize}
\item[($A_r$)] $V = \{v \in {\mathbb R}^{r+1} \mid \langle v,e_1 + \ldots +
e_{r+1} \rangle = 0\}$, $r \geq 1$;\\
$\Sigma = \{e_i - e_j \mid i \neq j\}$;
$\Sigma^+ = \{e_i - e_j \mid i < j\}$;\\
$\alpha_1 = e_1 - e_2 , \ldots,\alpha_r = e_r - e_{r+1}$;\\
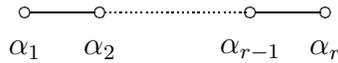
$$
\xy
\POS (0,0) *\cir<2pt>{} ="a",
(10,0) *\cir<2pt>{}="b",
(30,0) *\cir<2pt>{}="c",
(40,0) *\cir<2pt>{}="d",
(0,-5) *{\alpha_1},
(10,-5) *{\alpha_2},
(30,-5) *{\alpha_{r-1}},
(40,-5) *{\alpha_r},
\ar @{-} "a";"b",
\ar @{.} "b";"c",
\ar @{-} "c";"d",
\endxy
$$

\medskip\noindent
$M = SL_{r+1}({\mathbb R})/SO_{r+1}$: $(1,\ldots,1)$;\\
$M = SL_{r=1}({\mathbb C})/SU_{r+1}$: $(2,\ldots,2)$;\\
$M = SL_{r+1}({\mathbb H})/Sp_{r+1}$: $(4,\ldots,4)$;\\
$M = E_6^{-26}/F_4$: $(8,8)$;\\
$M = SO^o_{n+1,1}/SO_{n+1}$: $(n)$, $n \geq 2$;

\medskip
\item[($B_r$)] $V = {\mathbb R}^r$, $r \geq 2$;\\
$\Sigma = \{\pm e_i \pm e_j \mid i < j\} \cup \{\pm e_i\}$;
$\Sigma^+ = \{e_i \pm e_j \mid i < j\} \cup \{e_i\}$;\\
$\alpha_1 = e_1 - e_2 , \ldots,\alpha_{r-1} = e_{r-1} - e_r,\alpha_r=e_r$;\\
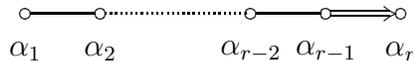
$$
\xy
\POS (0,0) *\cir<2pt>{} ="a",
(10,0) *\cir<2pt>{}="b",
(30,0) *\cir<2pt>{}="c",
(40,0) *\cir<2pt>{}="d",
(50,0) *\cir<2pt>{}="e",
(0,-5) *{\alpha_1},
(10,-5) *{\alpha_2},
(30,-5) *{\alpha_{r-2}},
(40,-5) *{\alpha_{r-1}},
(50,-5) *{\alpha_r},
\ar @{-} "a";"b",
\ar @{.} "b";"c",
\ar @{-} "c";"d",
\ar @2{->} "d";"e",
\endxy
$$

\medskip\noindent
$M = SO_{2r+1}({\mathbb C})/SO_{2r+1}$: $(2,\ldots,2,2)$;\\
$M = SO^o_{r+n,r}/SO_{r+n} SO_r$: $(1,\ldots,1,n)$, $n \geq 1$;

\medskip
\item[($C_r$)] $V = {\mathbb R}^r$, $r \geq 2$;\\
$\Sigma = \{\pm e_i \pm e_j \mid i < j\} \cup \{\pm 2e_i\}$;
$\Sigma^+ = \{e_i \pm e_j \mid i < j\} \cup \{2e_i\}$;\\
$\alpha_1 = e_1 - e_2 , \ldots,\alpha_{r-1} = e_{r-1} - e_r,\alpha_r=2e_r$;\\
$$
\xy
\POS (0,0) *\cir<2pt>{} ="a",
(10,0) *\cir<2pt>{}="b",
(30,0) *\cir<2pt>{}="c",
(40,0) *\cir<2pt>{}="d",
(50,0) *\cir<2pt>{}="e",
(0,-5) *{\alpha_1},
(10,-5) *{\alpha_2},
(30,-5) *{\alpha_{r-2}},
(40,-5) *{\alpha_{r-1}},
(50,-5) *{\alpha_r},
\ar @{-} "a";"b",
\ar @{.} "b";"c",
\ar @{-} "c";"d",
\ar @2{<-} "d";"e",
\endxy
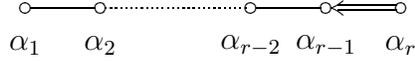
$$

\medskip\noindent
$M = Sp_r({\mathbb R})/U_r$: $(1,\ldots,1,1)$;\\
$M = Sp_r({\mathbb C})/Sp_r$: $(2,\ldots,2,2)$;\\
$M = Sp_{r,r}/Sp_r Sp_r$: $(4,\ldots,4,3)$;\\
$M = SU_{r,r}/S(U_r U_r)$: $(2,\ldots,2,1)$;\\
$M = SO_{2r}({\mathbb H})/U_{2r}$: $(4,\ldots,4,1)$;\\
$M = E_7^{-25}/E_6U_1$: $(8,8,1)$;

\medskip
\item[($D_r$)] $V = {\mathbb R}^r$, $r \geq 3$;\\
$\Sigma = \{\pm e_i \pm e_j \mid i < j\}$;
$\Sigma^+ = \{e_i \pm e_j \mid i < j\}$;\\
$\alpha_1 = e_1 - e_2 , \ldots,\alpha_{r-1} = e_{r-1} -
e_r,\alpha_r=e_{r-1}+e_r$;\\
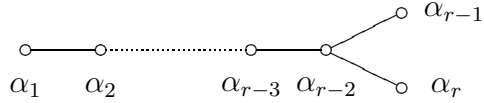
$$
\xy
\POS (0,0) *\cir<2pt>{} ="a",
(10,0) *\cir<2pt>{}="b",
(30,0) *\cir<2pt>{}="c",
(40,0) *\cir<2pt>{}="d",
(50,5) *\cir<2pt>{}="e",
(50,-5) *\cir<2pt>{}="f",
(0,-5) *{\alpha_1},
(10,-5) *{\alpha_2},
(30,-5) *{\alpha_{r-3}},
(40,-5) *{\alpha_{r-2}},
(57,5) *{\alpha_{r-1}},
(56,-5) *{\alpha_r},
\ar @{-} "a";"b",
\ar @{.} "b";"c",
\ar @{-} "c";"d",
\ar @{-} "d";"e",
\ar @{-} "d";"f",
\endxy
$$

\medskip\noindent
$M = SO^o_{r,r}/SO_r SO_r$: $(1,\ldots,1)$;\\
$M = SO_{2r}({\mathbb C})/SO_{2r}$: $(2,\ldots,2)$;

\medskip
\item[($E_6$)] $V = \{v \in {\mathbb R}^8 \mid \langle v,e_6-e_7 \rangle =
\langle v , e_7 + e_8 \rangle = 0\}$;\\
$\Sigma = \{\pm e_i \pm e_j \mid i < j \leq 5\} \cup \{ \frac{1}{2} \sum_{i=1}^8 (-1)^{n(i)}e_i \in V \mid \sum_{i=1}^8 (-1)^{n(i)}
\ {\rm even}\}$;\\
$\Sigma^+ = \{e_i \pm e_j \mid i > j\} \cup \{ \frac{1}{2} (e_8-e_7-e_6+\sum_{i=1}^5 (-1)^{n(i)}e_i) \mid \sum_{i=1}^5 (-1)^{n(i)}
\ {\rm even}\}$;\\
$\alpha_1 = \frac{1}{2}(e_1-e_2-e_3-e_4-e_5-e_6-e_7+e_8)$, $\alpha_2
= e_1 + e_2$,
$\alpha_i = e_{i-1} - e_{i-2} \ (i=3,\ldots,6)$;\\
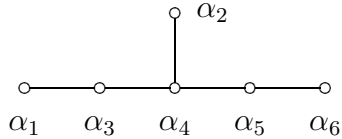
$$
\xy
\POS (0,0) *\cir<2pt>{} ="a",
(20,10) *\cir<2pt>{} = "b",
(10,0) *\cir<2pt>{}="c",
(20,0) *\cir<2pt>{}="d",
(30,0) *\cir<2pt>{}="e",
(40,0) *\cir<2pt>{}="f",
(0,-5) *{\alpha_1},
(25,10) *{\alpha_2},
(10,-5) *{\alpha_3},
(20,-5) *{\alpha_4},
(30,-5) *{\alpha_5},
(40,-5) *{\alpha_6},
\ar @{-} "a";"c",
\ar @{-} "c";"d",
\ar @{-} "b";"d",
\ar @{-} "d";"e",
\ar @{-} "e";"f",
\endxy
$$

\medskip\noindent
$M = E_6^6/Sp_4$: $(1,1,1,1,1,1)$;\\
$M = E_6^{\mathbb C}/E_6$: $(2,2,2,2,2,2)$;

\medskip
\item[($E_7$)] $V = \{v \in {\mathbb R}^8 \mid \langle v , e_7 + e_8 \rangle = 0\}$;\\
$\Sigma = \{\pm e_i \pm e_j \mid i < j \leq 6\} \cup \{ \pm
(e_7-e_8)\} \cup
\{ \frac{1}{2} \sum_{i=1}^8 (-1)^{n(i)}e_i \in V \mid \sum_{i=1}^8 (-1)^{n(i)}
\ {\rm even}\}$;\\
$\Sigma^+ = \{e_i \pm e_j \mid i > j\} \cup \{e_8 - e_7\} \cup
\{ \frac{1}{2}
(e_8-e_7+\sum_{i=1}^6 (-1)^{n(i)}e_i) \mid \sum_{i=1}^6 (-1)^{n(i)}
\ {\rm odd}\}$;\\
$\alpha_1 = \frac{1}{2}(e_1-e_2-e_3-e_4-e_5-e_6-e_7+e_8)$, $\alpha_2
= e_1 + e_2$,
$\alpha_i = e_{i-1} - e_{i-2} \ (i=3,\ldots,7)$;\\
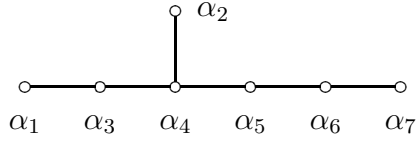
$$
\xy
\POS (0,0) *\cir<2pt>{} ="a",
(20,10) *\cir<2pt>{} = "b",
(10,0) *\cir<2pt>{}="c",
(20,0) *\cir<2pt>{}="d",
(30,0) *\cir<2pt>{}="e",
(40,0) *\cir<2pt>{}="f",
(50,0) *\cir<2pt>{}="g",
(0,-5) *{\alpha_1},
(25,10) *{\alpha_2},
(10,-5) *{\alpha_3},
(20,-5) *{\alpha_4},
(30,-5) *{\alpha_5},
(40,-5) *{\alpha_6},
(50,-5) *{\alpha_7},
\ar @{-} "a";"c",
\ar @{-} "c";"d",
\ar @{-} "b";"d",
\ar @{-} "d";"e",
\ar @{-} "e";"f",
\ar @{-} "f";"g",
\endxy
$$

\medskip\noindent
$M = E_7^7/SU_8$: $(1,1,1,1,1,1,1)$;\\
$M = E_7^{\mathbb C}/E_7$: $(2,2,2,2,2,2,2)$;

\medskip
\item[($E_8$)] $V = {\mathbb R}^8$;\\
$\Sigma = \{\pm e_i \pm e_j \mid i < j \} \cup \{\frac{1}{2}
\sum_{i=1}^8 (-1)^{n(i)}e_i \mid \sum_{i=1}^8 (-1)^{n(i)}
\ {\rm even}\}$;\\
$\Sigma^+ = \{e_i \pm e_j \mid i > j\} \cup
\{ \frac{1}{2}
(e_8+\sum_{i=1}^7 (-1)^{n(i)}e_i) \mid \sum_{i=1}^7 (-1)^{n(i)}
\ {\rm even}\}$;\\
$\alpha_1 = \frac{1}{2}(e_1-e_2-e_3-e_4-e_5-e_6-e_7+e_8)$, $\alpha_2
= e_1 + e_2$,
$\alpha_i = e_{i-1} - e_{i-2} \ (i=3,\ldots,8)$;\\
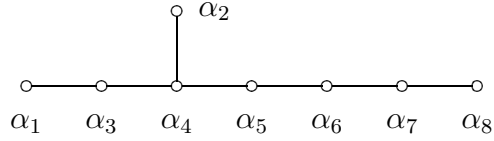
$$
\xy
\POS (0,0) *\cir<2pt>{} ="a",
(20,10) *\cir<2pt>{} = "b",
(10,0) *\cir<2pt>{}="c",
(20,0) *\cir<2pt>{}="d",
(30,0) *\cir<2pt>{}="e",
(40,0) *\cir<2pt>{}="f",
(50,0) *\cir<2pt>{}="g",
(60,0) *\cir<2pt>{}="h",
(0,-5) *{\alpha_1},
(25,10) *{\alpha_2},
(10,-5) *{\alpha_3},
(20,-5) *{\alpha_4},
(30,-5) *{\alpha_5},
(40,-5) *{\alpha_6},
(50,-5) *{\alpha_7},
(60,-5) *{\alpha_8},
\ar @{-} "a";"c",
\ar @{-} "c";"d",
\ar @{-} "b";"d",
\ar @{-} "d";"e",
\ar @{-} "e";"f",
\ar @{-} "f";"g",
\ar @{-} "g";"h",
\endxy
$$

\medskip\noindent
$M = E_8^8/SO_{16}$: $(1,1,1,1,1,1,1,1)$;\\
$M = E_8^{\mathbb C}/E_8$: $(2,2,2,2,2,2,2,2)$;

\medskip
\item[($F_4$)] $V = {\mathbb R}^4$;\\
$\Sigma = \{\pm e_i \pm e_j \mid i < j\} \cup \{ \pm e_i \} \cup
\{\frac{1}{2}(\pm e_1 \pm e_2 \pm e_3 \pm e_4\}$;\\
$\Sigma^+ = \{e_i \pm e_j \mid i < j\} \cup \{  e_i \} \cup
\{\frac{1}{2}(e_1 \pm e_2 \pm e_3 \pm e_4\}$;\\
$\alpha_1 = e_2 - e_3,\alpha_2 = e_3-e_4,\alpha_3=e_4,
\alpha_4 = \frac{1}{2}(e_1-e_2-e_3-e_4)$;\par
$$
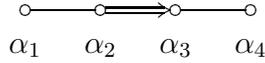
\xy
\POS (0,0) *\cir<2pt>{} ="a",
(10,0) *\cir<2pt>{}="b",
(20,0) *\cir<2pt>{}="c",
(30,0) *\cir<2pt>{}="d",
(0,-5) *{\alpha_1},
(10,-5) *{\alpha_2},
(20,-5) *{\alpha_3},
(30,-5) *{\alpha_4},
\ar @{-} "a";"b",
\ar @2{->} "b";"c",
\ar @{-} "c";"d",
\endxy
$$

\medskip\noindent
$M = F_4^4/Sp_3 Sp_1$: $(1,1,1,1)$;\\
$M = F_4^{\mathbb C}/F_4$: $(2,2,2,2)$;\\
$M = E_6^2/SU_6 Sp_1$: $(1,1,2,2)$;\\
$M = E_7^{-5}/SO_{12} Sp_1$: $(1,1,4,4)$;\\
$M = E_8^{-24}/E_7 Sp_1$: $(1,1,8,8)$;

\medskip
\item[($G_2$)] $V = \{v \in {\mathbb R}^3 \mid \langle v , e_1+e_2+e_3
\rangle = 0\}$;\\
$\Sigma = \{\pm (e_i - e_j) \mid i < j\} \cup
\{\pm(2e_i - e_j - e_k) \mid i \neq j \neq k \neq i\}$;\\
$\Sigma^+ =
\{e_1-e_2,-2e_1+e_2+e_3,-e_1+e_3,-e_2+e_3,-2e_2+e_1+e_3,2e_3-e_1-e_2\}$;\\
$\alpha_1 = e_1-e_2,\alpha_2 = -2e_1+e_2+e_3$;
$$
\xy
\POS (0,0) *\cir<2pt>{} ="a",
(10,0) *\cir<2pt>{}="b",
(0,-5) *{\alpha_1},
(10,-5) *{\alpha_2},
\ar @3{<-} "a";"b",
\endxy
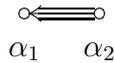
$$

\medskip\noindent
$M = G_2^2/SO_4$: $(1,1)$;\\
$M = G_2^{\mathbb C}/G_2$: $(2,2)$;

\medskip
\item[($BC_r$)] $V = {\mathbb R}^r$, $r \geq 1$;\\
$\Sigma = \{\pm e_i \pm e_j \mid i < j\} \cup \{\pm e_i\} \cup \{\pm 2e_i\}$;
$\Sigma^+ = \{e_i \pm e_j \mid i < j\} \cup \{e_i\} \cup \{2e_i\}$;\\
$\alpha_1 = e_1 - e_2 , \ldots,\alpha_{r-1} = e_{r-1} - e_r,\alpha_r=e_r$;
$$
\xy
\POS (0,0) *\cir<2pt>{} ="a",
(10,0) *\cir<2pt>{}="b",
(30,0) *\cir<2pt>{}="c",
(40,0) *\cir<2pt>{}="d",
(50,0) *\cir<2pt>{},
(50,0) *\cir<4pt>{}="e",
(0,-5) *{\alpha_1},
(10,-5) *{\alpha_2},
(30,-5) *{\alpha_{r-2}},
(40,-5) *{\alpha_{r-1}},
(53,-5) *{(\alpha_r,2\alpha_r)},
\ar @{-} "a";"b",
\ar @{.} "b";"c",
\ar @{-} "c";"d",
\ar @2{<->} "d";"e",
\endxy
$$

\medskip\noindent
$M = SU_{r+n,r}/S(U_{r+n} U_r)$: $(2,\ldots,2,(2n,1))$, $n
\geq 1$;\\
$M = SO_{2r+1}({\mathbb H})/U_{2r+1}$: $(4,\ldots,4,(4,1))$;\\
$M = Sp_{r+n,r}/Sp_{r+n} Sp_r$: $(4,\ldots,4,(4n,3))$, $n
\geq 1$;\\
$M = E_6^{-14}/Spin_{10} U_1$: $(6,(8,1))$;\\
$M = F_4^{-20}/Spin_9$: $((8,7))$.
\end{itemize}

\section{Parabolic subalgebras}

A subalgebra ${\mathfrak q}$ of ${\mathfrak g}$ is called parabolic
if there exists a maximal solvable subalgebra ${\mathfrak b}$ of ${\mathfrak g}$
such that ${\mathfrak b} \subset {\mathfrak q}$. A maximal solvable subalgebra
of ${\mathfrak g}$ is also known as a Borel subalgebra. The maximal solvable
subalgebras of real semisimple Lie algebras were classified by Mostow in \cite{M61}.
An example of a parabolic subalgebra is ${\mathfrak q} = {\mathfrak k}_0
\oplus {\mathfrak a} \oplus {\mathfrak n}$. This subalgebra is called minimal parabolic
as every parabolic subalgebra of ${\mathfrak g}$ contains a
subalgebra which is conjugate to ${\mathfrak k}_0
\oplus {\mathfrak a} \oplus {\mathfrak n}$.

We will now give a complete description of the parabolic subalgebras of ${\mathfrak g}$.
Let $\Phi$ be a
subset of $\Lambda$. We denote by $\Sigma_\Phi$ the root subsystem
of $\Sigma$ generated by $\Phi$, that is, $\Sigma_\Phi$ is the
intersection of $\Sigma$ and the linear span of $\Phi$, and put
$\Sigma_\Phi^+ = \Sigma_\Phi \cap \Sigma^+$.
We define an abelian subalgebra ${\mathfrak a}_\Phi$ of ${\mathfrak g}$ by
$$
{\mathfrak a}_\Phi = \bigcap_{\alpha \in \Phi} {\rm ker}\,\alpha. $$
Note that ${\mathfrak a}_\Phi = {\mathfrak a}$ if $\Phi = \emptyset$
and ${\mathfrak a}_\Phi = \{0\}$ if $\Phi = \Lambda$.
The normalizer $N_{\mathfrak g}({\mathfrak a}_\Phi)$ of ${\mathfrak a}_\Phi$
in ${\mathfrak g}$ is equal to the centralizer $Z_{\mathfrak g}({\mathfrak a}_\Phi)$
of ${\mathfrak a}_\Phi$ in ${\mathfrak g}$, and
$$
{\mathfrak l}_\Phi = N_{\mathfrak g}({\mathfrak a}_\Phi)
= Z_{\mathfrak g}({\mathfrak a}_\Phi) = {\mathfrak g}_0 \oplus \left(\bigoplus_{\lambda
\in \Sigma_\Phi} {\mathfrak g}_{\lambda}\right)
$$
is a reductive subalgebra of ${\mathfrak g}$. We define a nilpotent
subalgebra ${\mathfrak n}_\Phi$ of ${\mathfrak g}$ by
$$
{\mathfrak n}_\Phi = \bigoplus_{\lambda \in \Sigma^+\setminus
\Sigma_\Phi^+} {\mathfrak g}_{\lambda}.
$$
It is easy to verify from properties of root spaces that
$[{\mathfrak l}_\Phi,{\mathfrak n}_\Phi] \subset {\mathfrak n}_\Phi$,
and therefore
\[
{\mathfrak q}_\Phi = {\mathfrak l}_\Phi \oplus {\mathfrak n}_\Phi
\]
is a subalgebra of ${\mathfrak g}$, the so-called parabolic
subalgebra of ${\mathfrak g}$ associated with the subset $\Phi$ of
$\Lambda$. The decomposition ${\mathfrak q}_\Phi =
{\mathfrak l}_\Phi \oplus {\mathfrak n}_\Phi$ is known as the
Chevalley decomposition of the parabolic subalgebra ${\mathfrak
q}_\Phi$.

We denote by ${\mathfrak a}^\Phi = {\mathfrak a} \ominus {\mathfrak a}_\Phi$ the
orthogonal complement of ${\mathfrak a}_\Phi$ in ${\mathfrak a}$, and
define a reductive subalgebra ${\mathfrak m}_\Phi$ of
${\mathfrak g}$ by
\[
{\mathfrak m}_\Phi = {\mathfrak l}_\Phi \ominus {\mathfrak a}_\Phi =
{\mathfrak k}_0 \oplus {\mathfrak a}^\Phi \oplus
\left(\bigoplus_{\lambda \in \Sigma_\Phi} {\mathfrak
g}_{\lambda}\right)
\]
and a semisimple subalgebra ${\mathfrak g}_\Phi$ of ${\mathfrak m}_\Phi$ by
\[
{\mathfrak g}_\Phi = [{\mathfrak m}_\Phi,{\mathfrak m}_\Phi] =
[{\mathfrak l}_\Phi,{\mathfrak l}_\Phi].
\]
The center
${\mathfrak z}_\Phi$ of ${\mathfrak m}_\Phi$ is contained in
${\mathfrak k}_0$ and induces the direct sum decomposition
${\mathfrak m}_\Phi = {\mathfrak z}_\Phi \oplus {\mathfrak g}_\Phi$.
The decomposition
\[
{\mathfrak q}_\Phi = {\mathfrak m}_\Phi \oplus {\mathfrak a}_\Phi
\oplus {\mathfrak n}_\Phi
\]
is known as the Langlands decomposition of the parabolic subalgebra
${\mathfrak q}_\Phi$.

Each parabolic subalgebra of ${\mathfrak g}$ is conjugate in ${\mathfrak
g}$ to ${\mathfrak q}_\Phi$ for some subset $\Phi$ of $\Lambda$. The
set of conjugacy classes of parabolic subalgebras of ${\mathfrak g}$
therefore has $2^r$ elements. Two parabolic subalgebras ${\mathfrak
q}_{\Phi_1}$ and ${\mathfrak q}_{\Phi_2}$ of ${\mathfrak g}$ are
conjugate in the full automorphism group ${\rm Aut}({\mathfrak g})$
of ${\mathfrak g}$ if and only if there exists an automorphism $F$
of the Dynkin diagram associated to $\Lambda$ with $F(\Phi_1) =
\Phi_2$.

For each $\lambda \in \Sigma$ we define
\[
{\mathfrak k}_\lambda = {\mathfrak k} \cap ({\mathfrak g}_\lambda
\oplus {\mathfrak g}_{-\lambda})\ \ {\rm and}\ \ {\mathfrak
p}_\lambda = {\mathfrak p} \cap ({\mathfrak g}_\lambda \oplus
{\mathfrak g}_{-\lambda}).
\]
Then we have ${\mathfrak p}_\lambda = {\mathfrak p}_{-\lambda}$,
${\mathfrak k}_\lambda = {\mathfrak k}_{-\lambda}$ and ${\mathfrak
p}_\lambda \oplus {\mathfrak k}_\lambda = {\mathfrak g}_\lambda
\oplus {\mathfrak g}_{-\lambda}$ for all $\lambda \in \Sigma$. It is
easy to see that the subspaces
\[
{\mathfrak p}_\Phi = {\mathfrak l}_\Phi \cap {\mathfrak p} =
{\mathfrak a} \oplus \left( \bigoplus_{\lambda \in \Sigma_{\Phi}}
{\mathfrak p}_\lambda \right) \ {\rm and}\ {\mathfrak p}_\Phi^s =
{\mathfrak m}_\Phi \cap {\mathfrak p} = {\mathfrak g}_\Phi \cap
{\mathfrak p} = {\mathfrak a}^\Phi \oplus \left( \bigoplus_{\lambda
\in \Sigma_{\Phi}} {\mathfrak p}_\lambda \right)
\]
are Lie triple systems in ${\mathfrak p}$. We define a subalgebra
${\mathfrak k}_\Phi$ of ${\mathfrak k}$ by
\[
{\mathfrak k}_\Phi = {\mathfrak q}_\Phi \cap {\mathfrak k} =
{\mathfrak l}_\Phi \cap {\mathfrak k} = {\mathfrak m}_\Phi \cap
{\mathfrak k} = {\mathfrak k}_0 \oplus \left( \bigoplus_{\lambda \in
\Sigma_{\Phi}} {\mathfrak k}_\lambda \right).
\]
Then $${\mathfrak g}_\Phi = ({\mathfrak g}_\Phi \cap {\mathfrak
k}_\Phi) \oplus {\mathfrak p}_\Phi^s$$ is a Cartan decomposition of
the semisimple subalgebra ${\mathfrak g}_\Phi$ of ${\mathfrak g}$
and ${\mathfrak a}^\Phi$ is a maximal abelian subspace of
${\mathfrak p}_\Phi^s$. If we define $$({\mathfrak g}_\Phi)_0 =
({\mathfrak g}_\Phi \cap {\mathfrak k}_0) \oplus {\mathfrak
a}^\Phi,$$ then $${\mathfrak g}_\Phi = ({\mathfrak g}_\Phi)_0 \oplus
\left(\bigoplus_{\lambda \in \Sigma_\Phi} {\mathfrak
g}_{\lambda}\right)$$ is the restricted root space decomposition of
${\mathfrak g}_\Phi$ with respect to ${\mathfrak a}^\Phi$ and $\Phi$
is the corresponding set of simple roots. Since ${\mathfrak m}_\Phi
= {\mathfrak z}_\Phi \oplus {\mathfrak g}_\Phi$ and ${\mathfrak
z}_\Phi \subset {\mathfrak k}_0$, we see that ${\mathfrak g}_\Phi
\cap {\mathfrak k}_0 = {\mathfrak k}_0 \ominus {\mathfrak z}_\Phi$.

{\bf Example.} $M = SL_{r+1}({\mathbb R})/SO_{r+1}$.
For $\Phi = \emptyset$ we get the minimal parabolic subalgebra
$$
{\mathfrak q}_\emptyset = {\mathfrak a} \oplus {\mathfrak n},
$$
which is the solvable subalgebra of ${\mathfrak s}{\mathfrak l}_{r+1}({\mathbb R})$
consisting of all upper triangular $(r+1) \times (r+1)$-matrices with trace zero.

If $\phi = \Lambda$ we of course get
$$
{\mathfrak q}_\Lambda = {\mathfrak s}{\mathfrak l}_{r+1}({\mathbb R}).
$$

We now discuss the more interesting case of $\Phi \notin \{\emptyset,\Lambda\}$.
We denote by $r_\Phi$ the cardinality of the set $\Phi$. To begin with, we assume
that $\Phi$ is connected in the sense that either $r_\Phi = 1$, or $r_\Phi \geq 2$ and
$\alpha_{i-1} \in \Phi$ or $\alpha_{i+1} \in \Phi$ whenever $\alpha_i \in \Phi$.
Then the reductive subalgebra ${\mathfrak l}_\Phi$ is isomorphic to
${\mathfrak s}{\mathfrak l}_{r_\Phi+1}({\mathbb R}) \oplus {\mathbb R}^{r-r_\Phi}$.
The abelian subalgebra ${\mathfrak a}_\Phi$ coincides with the split component
${\mathbb R}^{r-r_\Phi}$ of ${\mathfrak l}_\Phi$, and the reductive subalgebra
${\mathfrak m}_\Phi$ corresponds to the semisimple subalgebra isomorphic to
${\mathfrak s}{\mathfrak l}_{r_\Phi+1}({\mathbb R})$, which in this case also coincides
with ${\mathfrak g}_\Phi$.

If $\Phi$
consists of several connected components, then ${\mathfrak l}_\Phi$ is isomorphic to
the direct sum of the special linear algebras of the corresponding
connected components of $\Phi$ and an abelian subalgebra of dimension $r-r_\Phi$.
In matrix form this direct sum
corresponds to the real vector space of $(r+1) \times (r+1)$-matrices with nonzero entries
only in a certain block diagonal decomposition. The nilpotent subalgebra
${\mathfrak n}_\Phi$ consists then of all strictly upper triangular matrices
all of whose entries in the intersection of the upper triangle and the blocks are zero.
Altogether this shows that the parabolic subalgebra ${\mathfrak q}_\Phi$
is given by matrices with zero entries below a certain block diagonal decomposition.
From this we conclude that the parabolic subalgebra
${\mathfrak q}_\Phi$ consists of all block diagonal upper triangular matrices,
where the block decomposition of the matrix corresponds to the decomposition of $\Phi$
into connected sets.

\section{Horospherical decompositions}

We now relate parabolic subalgebras to the geometry of the
symmetric space $M$. Let $\Phi$ be a subset of $\Lambda$ and $r_\Phi$
the cardinality of $\Phi$. We denote by $A_\Phi$ the connected abelian subgroup of
$G$ with Lie algebra ${\mathfrak a}_\Phi$ and by $N_\Phi$ the
connected nilpotent subgroup of $G$ with Lie algebra ${\mathfrak
n}_\Phi$. The centralizer $$L_\Phi = Z_G({\mathfrak a}_\Phi)$$ of
${\mathfrak a}_\Phi$ in $G$ is a reductive subgroup of $G$ with Lie
algebra ${\mathfrak l}_\Phi$. The subgroup $A_\Phi$ is contained in
the center of $L_\Phi$. The subgroup $L_\Phi$ normalizes $N_\Phi$
and $Q_\Phi = L_\Phi N_\Phi$ is a subgroup of $G$ with Lie algebra
${\mathfrak q}_\Phi$. The subgroup $Q_\Phi$ coincides with the
normalizer $N_G({\mathfrak l}_\Phi \oplus {\mathfrak n}_\Phi)$ of
${\mathfrak l}_\Phi \oplus {\mathfrak n}_\Phi$ in $G$, and hence
$Q_\Phi$ is a closed subgroup of $G$.
The subgroup
$$
Q_\Phi = L_\Phi N_\Phi = N_G({\mathfrak l}_\Phi \oplus {\mathfrak n}_\Phi)
$$
of $G$ is the parabolic subgroup of $G$ associated with the subsystem $\Phi$ of
$\Lambda$.

Let $G_\Phi$ be the connected subgroup of $G$ with Lie algebra
${\mathfrak g}_\Phi$. Since ${\mathfrak g}_\Phi$ is semisimple,
$G_\Phi$ is a semisimple subgroup of $G$. The intersection
$$K_\Phi = L_\Phi \cap K$$ is a maximal
compact subgroup of $L_\Phi$ and ${\mathfrak k}_\Phi$ is the Lie
algebra of $K_\Phi$. The adjoint group ${\rm Ad}(L_\Phi)$ normalizes
${\mathfrak g}_\Phi$, and consequently
$$M_\Phi = K_\Phi G_\Phi$$
is a
subgroup of $L_\Phi$. One can show that $M_\Phi$ is a closed
reductive subgroup of $L_\Phi$, $K_\Phi$ is a maximal compact
subgroup of $M_\Phi$, and the center $Z_\Phi$ of $M_\Phi$ is a
compact subgroup of $K_\Phi$. The Lie algebra of $M_\Phi$ is
${\mathfrak m}_\Phi$ and $L_\Phi$ is isomorphic to the Lie group
direct product $M_\Phi \times A_\Phi$, that is, $$L_\Phi = M_\Phi \times
A_\Phi.$$  The parabolic subgroup $Q_\Phi$ acts transitively on $M$
and the isotropy subgroup at $o$ is $K_\Phi$, that is, $$M =
Q_\Phi/K_\Phi.$$

Since ${\mathfrak g}_\Phi = ({\mathfrak g}_\Phi \cap {\mathfrak
k}_\Phi) \oplus {\mathfrak p}_\Phi^s$ is a Cartan decomposition of
the semisimple subalgebra ${\mathfrak g}_\Phi$, we have $[{\mathfrak
p}_\Phi^s,{\mathfrak p}_\Phi^s] = {\mathfrak g}_\Phi \cap {\mathfrak
k}_\Phi$. Thus $G_\Phi$ is the connected closed subgroup of $G$ with
Lie algebra $[{\mathfrak p}_\Phi^s,{\mathfrak p}_\Phi^s] \oplus
{\mathfrak p}_\Phi^s$. Since ${\mathfrak p}_\Phi^s$ is a Lie triple
system in ${\mathfrak p}$, the orbit $F_\Phi^s = G_\Phi \cdot o$ of
the $G_\Phi$-action on $M$ containing $o$ is a connected totally
geodesic submanifold of $M$ with $T_oF_\Phi^s = {\mathfrak
p}_\Phi^s$. If $\Phi = \emptyset$, then $F_\emptyset^s = \{o\}$,
otherwise $F_\Phi^s$ is a Riemannian symmetric space of noncompact
type and with ${\rm rank}(F_\Phi^s) = r_\Phi$, and
\[
F_\Phi^s = G_\Phi \cdot o = G_\Phi/(G_\Phi\cap K_\Phi) = M_\Phi
\cdot o = M_\Phi/K_\Phi.
\]
The submanifolds $F_\Phi^s$ are closely related to the
boundary components of $M$ in the context of the maximal Satake
compactification of $M$ (see e.g.\ \cite{BJ}).

Clearly, ${\mathfrak a}_\Phi$ is a Lie triple system as well, and
the corresponding totally geodesic submanifold is a Euclidean space
\[
{\mathbb E}^{r-r_\Phi} = A_\Phi \cdot o.
\]
Since the action of $A_\Phi$ on $M$ is free and $A_\Phi$ is simply
connected, we can identify ${\mathbb E}^{r-r_\Phi}$, $A_\Phi$ and
${\mathfrak a}_\Phi$ canonically.

Finally, ${\mathfrak p}_\Phi ={\mathfrak p}_\Phi^s \oplus {\mathfrak
a}_\Phi $ is a Lie triple system, and the corresponding totally
geodesic submanifold $F_\Phi$ is the symmetric space
\[
F_\Phi = L_\Phi \cdot o = L_\Phi/K_\Phi = (M_\Phi \times
A_\Phi)/K_\Phi =
 F_\Phi^s \times {\mathbb E}^{r-r_\Phi}.
\]

The submanifolds $F_\Phi$ and $F_\Phi^s$ have a nice geometric
interpretation. Choose $Z \in {\mathfrak a}$
such that $\alpha(Z) = 0$ for all $\alpha \in \Phi$ and $\alpha(Z) >
0$ for all $\alpha \in \Lambda \setminus \Phi$, and consider the
geodesic $\gamma_Z(t) = {\rm Exp}(tZ) \cdot o$ in $M$ with
$\gamma_Z(0) = o$ and $\dot{\gamma}_Z(0) = Z$. The totally geodesic
submanifold $F_\Phi$ is the union of all geodesics in $M$ which are
parallel to $\gamma_Z$, and $F_\Phi^s$ is the semisimple part of
$F_\Phi$ in the de Rham decomposition of $F_\Phi$ (see e.g.\
\cite{E96}, Proposition 2.11.4  and Proposition 2.20.10).

The group $Q_\Phi$ is diffeomorphic to the product $M_\Phi\times
A_\Phi\times N_\Phi$. This analytic diffeomorphism induces an
analytic diffeomorphism between $F_\Phi^s \times {\mathbb
E}^{r-r_\Phi} \times N_\Phi$ and $M$ known as a horospherical
decomposition of the symmetric space $M$:
$$
M \cong F_\Phi^s \times {\mathbb E}^{r-r_\Phi} \times N_\Phi.
$$
It turns out that these horospherical decompositions of symmetric
spaces of noncompact type provide a good framework for classifying
the hyperpolar homogeneous foliations.

{\bf Example.} $M = SL_{r+1}({\mathbb R})/SO_{r+1}$.
For $\Phi = \emptyset$ we obtain $F_\Phi = {\mathbb E}^r$, which is
a maximal flat in $M$. In particular we see that the semisimple part $F_\Phi^s$
of $F_\Phi$ consists of exactly one point. The subgroup $N_\Phi$ is a
maximal horocyclic subgroup of $SL_{r+1}({\mathbb R})$, and the orbits of
$N_\Phi$ give a foliation on $M$ by maximal horocycles. The horospherical decomposition in this
case is therefore given by the product of an $r$-dimensional Euclidean space
and a maximal horocycle.

For $\Phi = \Lambda$ we get $F_\Phi = M$, and therefore a trivial horospherical
decomposition.

The interesting case is whenever $\Phi$ is different from $\emptyset$ and $\Lambda$.
Assume that the cardinality of $\Phi$ is $r_\phi$, and $r_\Phi = r_\Phi^1 + \cdots + r_\Phi^k$,
where $k$ is the number of connected components of $\Phi$ and $r_\Phi^1,\ldots,r_\phi^k$ are
the cardinalities of the connected components. Then $F_\Phi^s$ is the Riemannian product of
$k$ smaller dimensional symmetric spaces of the same kind as $M$, that is,
$$
F_\Phi^s = SL_{r_\Phi^1 + 1}({\mathbb R})/SO_{r_\Phi^1+1} \times \cdots \times
SL_{r_\Phi^k + 1}({\mathbb R})/SO_{r_\Phi^k+1}.
$$
The horospherical decomposition of $M$ consists of the product of this $F_\Phi^s$, a
Euclidean space of dimension $r - r_\phi$, and a horocycle of a suitable dimension.

Of particular interest for us will be the case when the connected components of $\Phi$ are
all of cardinality one. In this case $F_\Phi^s$ is the Riemannian product of $k$ real hyperbolic
planes ${\mathbb R}H^2$.

\section{Homogeneous hyperpolar foliations on Euclidean spaces}

The totally geodesic subspaces of ${\mathbb E}^m$ are the affine subspaces
of the underlying vector space. Since affine subspaces are flat, ``polar''
and ``hyperpola'' have the same meaning in the Euclidean setting.

For each linear subspace $V$ of ${\mathbb E}^m$ we define a
foliation ${\mathcal F}_V^m$ on ${\mathbb E}^m$ by
\[
({\mathcal F_V^m})_p = p + V = \{ p + v \mid v \in V\}
\]
for all $p \in {\mathbb E}^m$. Geometrically, the leaves of the foliation
${\mathcal F}_V^m$ are the affine subspaces of
${\mathbb E}^m$ which are parallel to $V$. It is obvious that
${\mathcal F_V^m}$ is a hyperpolar homogeneous foliation on
${\mathbb E}^m$.

Every hyperpolar homogeneous foliation on
$\mathbb{E}^m$ is isometrically congruent to ${\mathcal F}_V^m$ for
some linear subspace $V$.

\section{Homogeneous hyperpolar foliations on hyperbolic spaces}

Let $M$ be
a hyperbolic space over a normed real division algebra ${\mathbb F} \in
\{{\mathbb R},{\mathbb C},{\mathbb H},{\mathbb O}\}$.
We denote such a hyperbolic space by ${\mathbb F} H^n$, where $n \geq 2$ is the dimension
of the manifold over the algebra ${\mathbb F}$, and $n = 2$ if ${\mathbb F} = {\mathbb O}$. We denote
by $G$ the connected component of the full isometry group of $M$. Then we have
\[
G =
\begin{cases}
SO^o_{n,1} & \text{if ${\mathbb F} = {\mathbb
R}$,} \\
SU_{n,1} & \text{if ${\mathbb F} = {\mathbb
C}$,} \\
Sp_{n,1} & \text{if ${\mathbb F} = {\mathbb
H}$,} \\
F_4^{-20} & \text{if ${\mathbb F} = {\mathbb O}$, $n = 2$}.
\end{cases}
\]
We denote by $K$ a maximal compact subgroup of $G$, which is unique up to conjugation,
and by $o \in M$ the unique fixed point of the action of $K$ on $M$.
Then ${\mathbb F} H^n$ can be identified with the homogeneous space $G/K$ in the usual way, and we
have
\[
{\mathbb F}H^n =
\begin{cases}
SO^o_{n,1}/SO_n & \text{if ${\mathbb F} = {\mathbb
R}$,} \\
SU_{n,1}/S(U_nU_1) & \text{if ${\mathbb F} = {\mathbb
C}$,} \\
Sp_{n,1}/Sp_nSp_1 & \text{if ${\mathbb F} = {\mathbb
H}$,} \\
F_4^{-20}/Spin_9 & \text{if ${\mathbb F} = {\mathbb O}$, $n = 2$}.
\end{cases}
\]
We denote by ${\mathfrak g}$ and ${\mathfrak k}$ the Lie algebra of $G$ and $K$ respectively,

As $M$ has rank one, there is exactly one simple root $\alpha$, and we have
\[
{\mathfrak g} = {\mathfrak g}_{-2\alpha} \oplus {\mathfrak g}_{-\alpha} \oplus
{\mathfrak g}_0 \oplus {\mathfrak g}_\alpha \oplus {\mathfrak g}_{2\alpha}.
\]
If ${\mathbb F} = {\mathbb R}$, then $\pm 2\alpha$ are not restricted roots and
hence ${\mathfrak g}_{\pm 2\alpha} = \{0\}$. The subalgebra ${\mathfrak n}
= {\mathfrak g}_\alpha \oplus {\mathfrak g}_{2\alpha}$ of ${\mathfrak g}$ is
nilpotent and ${\mathfrak a} \oplus {\mathfrak n}$ is a solvable
subalgebra of ${\mathfrak g}$. The vector space decomposition
${\mathfrak g} = {\mathfrak k} \oplus {\mathfrak a} \oplus
{\mathfrak n}$ is an Iwasawa decomposition of ${\mathfrak g}$. Since $M$
is isometric to the solvable Lie group $AN$ equipped with a
suitable left-invariant Riemannian metric, it is
obvious that every subalgebra of ${\mathfrak a} \oplus {\mathfrak
n}$ of codimension one induces a homogeneous codimension one
foliation on $M$. Here are two examples:

1. The subspace ${\mathfrak n}$ is a subalgebra of ${\mathfrak a}
\oplus {\mathfrak n}$ of codimension one. Therefore the orbits of
the action of the nilpotent Lie group $N$ on $M$ form a homogeneous
foliation of codimension one on $M$. This foliation is the
well-known horosphere foliation.

2. Let $\ell$ be a one-dimensional linear subspace of ${\mathfrak
g}_\alpha$. Then ${\mathfrak a} \oplus ({\mathfrak n} \ominus \ell)$ is a
subalgebra of ${\mathfrak a} \oplus {\mathfrak n}$ of codimension
one. The orbits of the corresponding connected closed subgroup of
$AN$ form a homogeneous foliation on $M$. Different choices of $\ell$
lead to isometrically congruent foliations. We denote by
${\mathcal F}_{\mathbb F}^n$ a representative of this foliation.
The leaves of ${\mathcal F}_{\mathbb R}^n$ are a totally geodesic hyperplane ${\mathbb
R}H^{n-1}$ in ${\mathbb R}H^n$ and the equidistant hypersurfaces. If ${\mathbb F} =
{\mathbb C}$, one of the laaves of ${\mathcal F}_{\mathbb C}^n$ is the minimal ruled real
hypersurface in ${\mathbb C}H^n$ generated by a horocycle in a
real hyperbolic plane embedded totally geodesically
into ${\mathbb C}H^n$ as a real surface; the other leaves are the equidistant hypersurfaces.

It was proved by the author and Hiroshi Tamaru in \cite{BT03} that on ${\mathbb
F}H^n$ every homogeneous foliation of codimension one is
isometrically congruent to one of these two foliations. Since the
rank of ${\mathbb F}H^n$ is equal to one, every hyperpolar foliation
must have codimension one, and therefore every homogeneous
hyperpolar foliation on ${\mathbb F}H^n$ is isometrically
congruent to either the horosphere foliation on ${\mathbb F}H^n$ or
to ${\mathcal F}_{\mathbb F}^n$.

\section{Homogeneous hyperpolar foliations on products of hyperbolic spaces}

Let
\[
M = {\mathbb F}_1H^{n_1} \times \ldots \times {\mathbb
F}_kH^{n_k}
\]
be the Riemannian product of $k$ hyperbolic spaces, where $k \geq 2$ is a positive integer
and ${\mathbb F}_k \in \{{\mathbb R},{\mathbb C},{\mathbb H},{\mathbb O}\}$. Then
\[
{\mathcal F}_{{\mathbb F}_1}^{n_1} \times \ldots \times {\mathcal
F}_{{\mathbb F}_k}^{n_k}
\]
is a hyperpolar homogeneous foliation on $M$. This is an elementary
consequence of the previous example.

\section{Hyperpolar homogeneous foliations on products
of hyperbolic
spaces and a Euclidean space}

Let
\[
M = {\mathbb F}_1H^{n_1} \times \ldots \times {\mathbb
F}_kH^{n_k} \times {\mathbb E}^{m}
\]
be the Riemannian product of
$k$ hyperbolic spaces and an $m$-dimensional
Euclidean space, where $k$ and $m$ are positive integers. Moreover,
let $V$ be a linear subspace of ${\mathbb E}^m$. Then
\[
{\mathcal F}_{{\mathbb F}_1}^{n_1} \times \ldots \times {\mathcal
F}_{{\mathbb F}_k}^{n_k} \times {\mathcal F}_V^m
\]
is a hyperpolar homogeneous foliation on $M$.

\section{Homogeneous foliations on symmetric spaces of noncompact type}

Let $M$ be a Riemannian
symmetric space of noncompact type and $\Phi$ be a subset of
$\Lambda$ with the property that any two roots in $\Phi$ are not
connected in the Dynkin diagram of the restricted root system
associated with $M$. We call such a subset $\Phi$ an
orthogonal subset of $\Lambda$.

Each simple root $\alpha \in \Phi$
determines a totally geodesic hyperbolic space ${\mathbb F}_\alpha
H^{n_\alpha} \subset M$. In fact, ${\mathbb F}_\alpha H^{n_\alpha}
\subset M$ is the orbit of $G_\Phi$, $\Phi = \{\alpha\}$, through the point $o$.
If $2\alpha \notin \Sigma$, that is, if the vertex in the Dynkin
diagram corresponding to $\alpha$ is of the form $ \xy  *\cir<2pt>{} \endxy $,
then ${\mathbb F}_\alpha = {\mathbb R}$ and the dimension $n_\alpha$ is equal
to $m_\alpha + 1$, where $m_\alpha$ is the multiplicity
of the root $\alpha$.
If $2\alpha \in \Sigma$, that is, if the vertex in the Dynkin
diagram corresponding to $\alpha$ is of the form $ \xy  *\cir<2pt>{} *\cir<4pt>{} \endxy $,
then ${\mathbb F}_\alpha \in \{{\mathbb C},{\mathbb H},{\mathbb O}\}$.
Note that this can happen only if $\Sigma$ is of type $(BC_r)$
and $\alpha = \alpha_r$. We have $m_{2\alpha} \in \{1,3,7\}$,
\[
{\mathbb F}_\alpha =
\begin{cases}
{\mathbb C} & \text{if $m_{2\alpha} = 1$,} \\
{\mathbb H} & \text{if $m_{2\alpha} = 3$,} \\
{\mathbb O} & \text{if $m_{2\alpha} = 7$,}
\end{cases}
\]
and
\[
n_\alpha =
\begin{cases}
\frac{m_\alpha}{2} + 1 & \text{if $m_{2\alpha} = 1$,} \\
\frac{m_\alpha}{4} + 1 & \text{if $m_{2\alpha} = 3$,} \\
2 & \text{if $m_{2\alpha} = 7$.}
\end{cases}
\]

The symmetric space $F_\Phi$ in the horospherical decomposition of $M$ induced from $\Phi$
is isometric to
the Riemannian product of $r_\Phi$ hyperbolic spaces
and an $(r-r_\Phi)$-dimensional Euclidean space, that is,
\[
F_\Phi = F_\Phi^s \times {\mathbb E}^{r-r_\Phi} \cong \left(
\prod_{\alpha \in \Phi} {\mathbb F}_{\alpha} H^{n_\alpha} \right)
\times {\mathbb E}^{r-r_\Phi}.
\]
Then
\[
{\mathcal F}_\Phi = \prod_{\alpha \in \Phi} {\mathcal F}_{{\mathbb
F}_\alpha}^{n_\alpha}.
\]
is a hyperpolar homogeneous foliation on ${\mathcal F}_\Phi^s$. Let
$V$ be a linear subspace of ${\mathbb E}^{r-r_\Phi}$. Then
\[
{\mathcal F}_{\Phi,V} = {\mathcal F}_\Phi \times {\mathcal
F}_V^{r-r_\Phi} \times N_\Phi \subset F_\Phi^s \times {\mathbb
E}^{r-r_\Phi} \times N_\Phi = F_\Phi \times N_\Phi \cong M
\]
is a homogeneous foliation on $M$.

Recall that each foliation ${\mathcal F}_{{\mathbb
F}_\alpha}^{n_\alpha}$ on ${\mathbb F}_{\alpha} H^{n_\alpha}$
corresponds to a subalgebra of ${\mathfrak g}_{\{\alpha\}}$ of the
form ${\mathfrak a}^{\{\alpha\}} \oplus ({\mathfrak g}_\alpha \ominus \ell_\alpha)
\oplus {\mathfrak g}_{2\alpha}$ with some one-dimensional linear subspace
$\ell_\alpha$ of ${\mathfrak g}_\alpha$. Thus the foliation
${\mathcal F}_\Phi$ on $F_\Phi^s$ corresponds to the subalgebra
$${\mathfrak a}^\Phi \oplus \left( \bigoplus_{\alpha \in \Phi} \left(
({\mathfrak g}_\alpha \ominus \ell_\alpha)  \oplus {\mathfrak
g}_{2\alpha} \right) \right) = {\mathfrak a}^\Phi \oplus ({\mathfrak
n}_\Phi \ominus \ell_\Phi)$$ of ${\mathfrak g}_\Phi$, where
$\ell_\Phi=\bigoplus_{\alpha\in\Phi}\ell_\alpha$. Therefore the
foliation ${\mathcal F}_{\Phi,V}$ on $M$ corresponds to the
subalgebra
\[
{\mathfrak s}_{\Phi,V}=({\mathfrak a}^\Phi \oplus V)
\oplus({\mathfrak n}_\Phi\ominus\ell_\Phi)=({\mathfrak a}^\Phi \oplus V
\oplus {\mathfrak n}_\Phi)\ominus\ell_\Phi \subset {\mathfrak a} \oplus
{\mathfrak n}_\Phi
\]
of ${\mathfrak q}_\Phi$, where we identify canonically $V \subset
{\mathbb E}^{r-r_\Phi} = A_\Phi \cdot o$ with the corresponding
subspace of ${\mathfrak a}_\Phi$.
One can show that different choices of $\ell_\alpha$
in ${\mathfrak g}_\alpha$ lead to isometrically congruent foliations
${\mathcal F}_{\Phi,V}$ on $M$.

\section{The Classification}

We are now in a position to formulate the main classification result.

\begin{theorem}\label{maintheorem}
Let $M$ be a connected Riemannian symmetric space of noncompact
type.
\begin{itemize}
\item[(i)] Let $\Phi$ be an orthogonal subset of $\Lambda$ and $V$ be
a linear subspace of ${\mathbb E}^{r-r_\Phi}$. Then
\[
{\mathcal F}_{\Phi,V} = {\mathcal F}_\Phi \times {\mathcal
F}_V^{r-r_\Phi} \times N_\Phi \subset F_\Phi^s \times {\mathbb
E}^{r-r_\Phi} \times N_\Phi = M
\]
is a hyperpolar homogeneous foliation on $M$. \item[(ii)] Every
hyperpolar homogeneous foliation on $M$ is isometrically congruent
to ${\mathcal F}_{\Phi,V}$ for some orthogonal subset $\Phi$ of
$\Lambda$ and some linear subspace $V$ of ${\mathbb E}^{r-r_\Phi}$.
\end{itemize}
\end{theorem}

For the proof we refer to \cite{BDT}. The special case for codimension one
foliations was already solved in \cite{BT03}.

{\bf Example.} $M = SL_{r+1}({\mathbb R})/SO_{r+1}$.
The Dynkin diagram associated with $M$ is
$$
\xy
\POS (0,0) *\cir<2pt>{} ="a",
(10,0) *\cir<2pt>{}="b",
(30,0) *\cir<2pt>{}="c",
(40,0) *\cir<2pt>{}="d",
(0,-5) *{\alpha_1},
(10,-5) *{\alpha_2},
(30,-5) *{\alpha_{r-1}},
(40,-5) *{\alpha_r},
\ar @{-} "a";"b",
\ar @{.} "b";"c",
\ar @{-} "c";"d",
\endxy
$$
and therefore
the orthogonal subsets $\Phi$ of $\Lambda$ correspond precisely to the subsets
of $\{1,\ldots,r\}$ containing no two adjacent positive integers. Let $k= r_\Phi$ be the cardinality
of an orthogonal subset $\Phi$ of $\Lambda$. Since the multiplicity
of each simple root is one, the semisimple part  $F_\Phi^s$ of the horospherical decomposition
$M \cong F_\Phi^s \times {\mathbb E}^{r-k} \times N_\Phi$
is isometric
to the Riemannian product of $k$ real hyperbolic planes ${\mathbb R}H^2$. On
each of these real hyperbolic planes we choose the foliation determined by a geodesic
and its equidistant curves. The product of these foliations determines the foliation ${\mathcal F}_\Phi$
on the $k$-fold product $F_\Phi^s$ of real hyperbolic planes. On the abelian part ${\mathbb E}^{r-k}$
we choose a foliation ${\mathcal
F}_V^{r-k}$ by parallel affine subspaces (including the trivial foliations of
dimension $0$ and $r-k$). The product foliation ${\mathcal F}_\Phi \times {\mathcal
F}_V^{r-k}$ on the totally geodesic submanifold $F_\Phi = F_\Phi^s \times {\mathbb E}^{r-k}$
of $M$ is hyperpolar. The foliation $F_{\Phi,V}$ is then obtained by taking the product of
this foliation with the horocycle foliation $N_\Phi$ on $M$. Theorem \ref{maintheorem} says
that every hyperpolar foliation on $M = SL_{r+1}({\mathbb R})/SO_{r+1}$ is obtained in this way.

\bigskip\bigskip\noindent
{\sc Author's Address:} \\ Department of Mathematics, King's College London, \\
Strand, London, WC2R 2LS, United Kingdom \\
{\sc Email:}\\ jurgen.berndt@kcl.ac.uk

\end{document}